\documentclass{amsart}
\usepackage{graphicx}
\usepackage{amssymb,latexsym,amsmath,amscd,epsfig,amsthm}
\newcommand {\PP}{{I\kern-.3em P}}
\newcommand {\ZZ}{{Z\kern-.45em Z}}

\newcommand {\RR}{{I\kern-.3em R}}
\newcommand {\R}{\mathbb{R}}
\newcommand {\CC}{{I\kern-.6em C}}
\newcommand {\C}{\mathbb{C}}
\newcommand {\D}{\mathbb{D}}
\newcommand {\eps}{\epsilon}

\newcommand{\zl}{\zeta^{(l)}}
\newcommand{\ql}{Q^{(l)}}
\newcommand{\ml}{\mu^{(l)}}
\newcommand{\dl}{\delta^{(l)}}
\newcommand{\cal}{\mathcal}
\newcommand{\ldist}{\mbox{dist}(z,\ql)}
\newcommand{\beq}{\begin{equation}}
\newcommand{\eeq}{\end{equation}}
\newtheorem{lemma}{Lemma}[section]
\newtheorem{theorem}[lemma]{Theorem}

\newtheorem{proposition}{Proposition}

\newtheorem{corollary}{Corollary}
\newtheorem*{corollary*}{Corollary}

\begin{document}

\title[Hilbert Lemniscate]{A Hilbert Lemniscate Theorem in $\C^2$ }%
\author{T. Bloom, N. Levenberg and Yu. Lyubarskii}%
\address{T.B.: University of Toronto, Toronto, CANADA}
\email{bloom@math.toronto.edu}
\address{N.L.: Indiana University, Bloomington, IN 47405 USA}
\email{nlevenbe@indiana.edu}
\address{Yu.L.: Norwegian Univ. of Science and Technology, Trondheim, NORWAY}
\email{yura@math.ntnu.no}
\thanks{Supported in part by an NSERC grant (TB) and by the Norwegian Research Council project 160192/V30 (YuL)  }%
\subjclass{32U05, 32W20}%
\keywords{Logarithmic potential, Monge-Amp\`ere measure, subharmonic functions, atomization}%





\maketitle


\section{\bf Introduction.}
    \label{sec:1}

Let $K\subset \C$ be a compact set with connected complement. The
{\it Hilbert lemniscate theorem} in one variable says that for such
sets, given any $\epsilon>0$, there exists a polynomial $p$ with
\beq\label{lem} K\subset {\cal K}_p:=\{z:|p(z)| \leq
||p||_K:=\sup_{z\in K}|p(z)|\} \subset K^{\epsilon} :=\{z:
\hbox{dist}(z,K) \leq \epsilon\}. \eeq The set ${\cal K}_p$ is
called a {\it lemniscate}. In general, given $\epsilon>0$, one can
take $p$ to be a {\it Fekete polynomial} of sufficiently large
degree. A Fekete polynomial of degree $n$ for $K$ is a monic
polynomial $F_n(z)=\prod_{j=1}^n(z-a_{nj})$ with $a_{nj}\in K$
chosen so that $$\prod_{j<k}^n|a_{nj}-a_{nk}|=\max_{z_1,...,z_n\in
K} \prod_{j<k}^n|z_j-z_k|.$$

The condition that $K$ have connected complement is equivalent to
the {\it polynomial convexity} of $K$: this means that $K=\hat K$
where
$$\hat K:=\{z\in \C: |p(z)|\leq ||p||_K:=\sup_{\zeta \in K}|p(\zeta)| \ \hbox{for all
polynomials} \ p\}.$$ (Here and in the entire paper ``polynomial''
means {\it holomorphic} polynomial). We call $K$ {\it regular} if
the extremal function \beq\label{vk} V_K(z):=\max [0,\sup
\{{1\over{\rm deg}p}\log {|p(z)|}: p \ \hbox{polynomial}, \ {\rm
deg}p\geq 1, \ ||p||_K\leq 1\}] \eeq is continuous on $\C$. For the
lemniscate ${\cal K}_p$ in (\ref{lem}),
$$V_{{\cal K}}(z)=\max[{1\over {\rm deg}p}\log {[|p(z)|/||p||_K]},0].$$ If $K$ is regular, in choosing, e.g., a sequence of Fekete polynomials $\{F_n\}$, the functions
\beq\label{fek} {1\over n}\log {[|F_n(z)|/||F_n||_K]}\to V_K(z) \eeq
locally uniformly outside of $K$. We also have the normalized
counting measure of the zeros \beq\label{muk} \mu_n:={1\over n}
\sum_{j=1}^n \delta_{a_{nj}}\to {1\over 2\pi} \Delta V_K \eeq weak-*
as measures. Here, $\Delta V_K$, the Laplacian of $V_K$, is to be
interpreted as a positive distribution, i.e., a positive measure.
Another example of a sequence of polynomials for which (\ref{fek})
and (\ref{muk}) hold is gotten by taking the interval $K=[-1,1]$ and
the classical Chebyshev polynomials $\{T_n\}$. Here $T_n(x)=\cos
{n(\arccos x)}$ for $x\in \R$; $V_K(z)=\log |z+\sqrt {z^2-1}|$ and
the normalized counting measure of the zeros approximate the arcsine
distribution ${dx\over \sqrt {1-x^2}}=\Delta V_K$.

In several complex variables, given a compact set $K\subset \C^N$,
$N>1$, we can define the extremal function $V_K$ as in (\ref{vk})
where $p(z)=p(z_1,...,z_N)$ is a polynomial of the complex variables
$z_1,...,z_N$. The definitions of regularity and polynomial
convexity are defined as in the one-variable case; however this
latter definition is no longer equivalent to the complement of $K$
being connected. It follows from the definition of $V_K$ and $\hat
K$ that $V_K=V_{\hat K}$ and that $\hat K=\{z: V_K(z)=0\}$ so that
an assumption of polynomial convexity is a natural one. In this
paper, we will prove a version of Hilbert's lemniscate theorem in
$\C^2$, including a convergence of measures result in the spirit of
(\ref{muk}).

To motivate this result, we note that in several complex variables,
sublevel sets $\{z: |p(z)| \leq M\}$ for a polynomial $p$ are
unbounded; in general, one needs at least $N$ polynomials
$p_1,...,p_N$ to have hopes of a sublevel set $\{z\in \C^N: |p_1(z)|
\leq M_1,...,|p_N(z)| \leq M_N\}$ being compact. Moreover, the
topology of such sublevel sets can be complicated. A {\it polynomial
polyhedron} is a set $P$ which is the closure of the union of a
finite number of connected components of
$${\cal P}:=\{z\in \C^N: |p_1(z)| < 1, ... , |p_m(z)|<  1\}$$
where $p_1,...,p_m$ are polynomials. It is an easy exercise to see
that {\it given any polynomially convex compact set $K\subset \C^N$,
and any open neighborhood $\Omega$ of $K$, there exists a set of the
form ${\cal P}$ with $K\subset {\cal P} \subset \Omega$} (cf.
\cite{H1}). What is not at all obvious is a deep result of Bishop
\cite{bishop}: {\it there exists a {\bf special} polynomial
polyhedron $P$ with the same property}. We call a polynomial
polyhedron $P\subset \C^N$ {\it special} if it can be defined by
{\it exactly} $N$ polynomials. We emphasize that not all components
of ${\cal P}$ need be included in $P$. It is known (cf. \cite{K},
Theorem 5.3.1) that if the set
$${\cal P}:=\{z\in \C^N: |p_1(z)|< 1,..., |p_N(z)|< 1\},$$
consisting of the union of {\it all} components of a special
polynomial polyhedron defined by $p_1,...,p_N$ with
deg$p_1=...$deg$p_N=:n$ is compact, and if $(p_1,...,p_N):\C^N \to
\C^N$ is proper, then we have
$$V_{\cal P}(z)=\max \bigl [{1\over n} \log {|p_1(z)|}, ..., {1\over n} \log {|p_N(z)|},0 \bigr ].$$
Thus, it will be helpful to know when a compact set $K$ can be
approximated not just by a special polynomial polyhedron $P$, but by
the full component set ${\cal P}$ of such an object. It turns out
that if we work in $\C^2$ with variables $(z,w)$ and we assume, in
addition to $K=\hat K$, that $K\subset \C^2$ is {\it circled}; i.e.,
$z\in K$ if and only if $e^{it}z\in K$, then such an approximation
is possible. Moreover, in this case, utilizing one-variable
techniques, we can construct Bishop-type approximants which satisfy
an analogue of (\ref{fek}) and (\ref{muk}).

\begin{theorem}
    \label{th:1.1} Let $K\subset \C^2$ be a regular, circled, polynomially convex compact set. Then there exists a sequence of pairs of homogeneous polynomials $\{P_n,Q_n\}$, deg$P_n=$deg$Q_n=n$ with no common linear factors such that
$$\tilde u_n(z,w):= \max [{1\over n} \log {|P_n(z,w)|}, {1\over n} \log {|Q_n(z,w)|},0]$$
uniformly approximates $V_K$ on $\C^2$;
$$U_n(z,w):= \max [{1\over n} \log {|P_n(z,w)-1|}, {1\over n} \log {|Q_n(z,w)-1|}]$$
locally uniformly approximates $V_K$ on $\C^2\setminus \partial K$;
and
$$(dd^c\tilde u_n)^2 \to (dd^cV_K)^2$$
weak-* as measures in $\C^2$. Moreover, if $K$ is the closure of a
strictly pseudoconvex domain (e.g., a ball), then
$$(dd^c U_n)^2 \to (dd^cV_K)^2.$$
\end{theorem}

Here, for certain plurisubharmonic (psh) functions $u$ in $\C^2$,
the complex Monge-Amp\`ere measure $(dd^cu)^2$ associated to $u$ is
well-defined. We discuss this issue in section 4. In particular, for
regular compact sets $K\subset \C^2$, $(dd^cV_K)^2$ plays a role
analogous to $\Delta V_K$ in one variable. In Theorem 1.1,
\begin{itemize}
\item the function $\tilde u_n$ is the extremal function for the set
\beq\label{calk} {\cal K}_n:=\{(z,w)\in \C^2: |P_n(z,w)|\leq 1, \
|Q_n(z,w)|\leq 1\}; \eeq
\item the Monge-Amp\`ere measure $(dd^cU_n)^2$ is supported on the finite point set (see section 4) \beq\label{kd} K_n:=\{(z,w):P_n(z,w)=Q_n(z,w)=1\};
\eeq
\item  the measures $\{(dd^c\tilde u_n)^2\}_{n=1,...}, \ \{(dd^cU_n)^2\}_{n=1,...}$ are supported in a fixed compact set in $\C^2$.
\end{itemize}

\noindent The distinction between the sequences $\{\tilde u_n\}$ and
$\{U_n\}$ can easily be seen even in one variable: take
$K=\D:=\{t\in \C: |t|\leq 1\}$, the closed unit disk. Then
$V_{\D}(t)=\max [\log |t|,0]$ and, taking $p_n(t)=t^n$, we have
$$\tilde v_n(t):=\max[{1\over n}\log |p_n(t)|,0]\equiv V_{\D}(t)$$
while
$$V_n(t):={1\over n}\log |p_n(t)-1|={1\over n}\log |t^n-1|$$
converges locally uniformly to $V_{\D}$ in $\C \setminus \{|t|=1\}$
but we clearly do not have $V_n \to V_{\D}$ pointwise, or even ``in
capacity'' (cf., \cite{stahl}) on the circle $\{|t|=1\}$. However,
we do have $V_n \to V_{\D}$ in $L^1_{loc}(\C)$. Thus, we can utilize
elementary distribution theory to conclude that the normalized
counting measure of the zeros of these Fekete polynomials $p_n(t)$
converge weak-* to $\Delta V_{\D}$. Of course, in this example, the
convergence of these measures is trivial (and, as mentioned earlier,
always holds for Fekete polynomials). We discuss the analogous
example of the unit bidisk in $\C^2$ in section 4.

We prove the first part of Theorem 1.1 by reducing it to a
one-variable approximation problem in  section 2. Given a measure
$\mu$ in  $\C$  with $\mu(\C)=1$ consider its logarithmic potential

 \beq
 V(t) = \int_{\C} \log |1- \frac t \zeta|
d\mu(\zeta).
 \eeq
We assume that
 \beq
\lim_{|t|\to \infty} [V(t) - \log |t|] \ \hbox{exists}, \eeq

\beq \int_{\C} |\log |t|| d\mu(t) < \infty, \eeq \noindent and that
$V(t)$ is continuous in $\C$. Under these assumptions, we will prove
the following theorem, which is of interest in its own right, in
section 3:

\begin{theorem}
    \label{th:1.2} Given $V$ satisfying (7), (8) and (9), for each $\eps>0$ there exist  a number $N$ and polynomials $p(t)$
and $q(t)$ of degree $N$ such that \beq |V(t) - \frac 1 N \max \{
\log|p(t)|, \log|q(t)|\}  |<\eps, \ \ \ t\in \C. \label{max_claim}
\eeq
\end{theorem}

\noindent The construction is based on techniques developed in
\cite{lyumal}. There the authors construct an $L^1-$approximant  to
an arbitrary subharmonic function $u$ in $\C$ of the form $\log |f|$
with a (single) entire function $f$. The proof utilizes a clever
partition of $\C$ related to the measure $\mu$ and its support, due
to Yulmukhametov \cite{yul}. The precise version of the result that
we use in section 2 is labeled Lemma A. We remark that the genesis
of Theorem 1.2 occurred during an Oberwolfach meeting attended by
the second and third authors in February 2004.

In the final section of the paper, we turn to the proof of
Monge-Amp\`ere convergence, the second part of Theorem 1.1. For the
sequence $\{\tilde u_n\}$ this convergence is automatic; but for the
sequence $\{U_n\}$, which is {\it not} locally bounded, a
non-trivial argument is required. This is given as Theorem 4.1. We
would like to thank Urban Cegrell for pointing out an error in our
proof of this result in a previous version.

We remark that from Bishop's theorem one can construct sequences of
psh functions with the same properties as the sequence $\{U_n\}$ in
Theorem 1.1 for general regular, polynomially convex compact sets
$K\subset {\bf C}^N$ which are not necessarily circled. However,
this work of Bishop is technically complicated and the construction
may not yield psh functions which are the maximum of exactly $N$
functions of the form $c\log |p|$ where $p$ is a polynomial. Our
methods in constructing the polynomials in Theorem 1.1 are purely
one-variable in nature and provide, via the sets $\{K_n\}$ in
(\ref{kd}), discrete approximations to the Monge-Amp\`ere measure
$(dd^cV_K)^2$.

\vskip 5mm

\section{\bf Reduction to one-variable.}
    \label{sec:2}

\vskip 3mm

For $N=1,2,...$, let
$$L(\C^N):=\{u \ \hbox{psh in} \ \C^N: \ u(z) \leq \log^+ |z| +C \}$$
denote the class of psh functions of logarithmic growth on $\C^N$
where the constant $C$ can depend on $u$. For example, given a
polynomial $p$, $u(z):={1\over{\rm deg}p}\log {|p(z)|}\in L(\C^N)$.
We also consider the class
$$L^+(\C^N):=\{u\in L(\C^N): \log^+ |z| +C_1 \leq u(z) \leq \log^+ |z|
+C_2, \ \hbox{some} \ C_1,C_2 \}.$$ Note functions in this class are
locally bounded.

For a bounded Borel set $E$ in $\C^N$, one can define
\beq\label{ve}V_E(z):=\sup \{u(z):u \in L(\C^N), \ u\leq 0 \
\hbox{on} \
    E\}.
    \eeq
The uppersemicontinuous (usc) regularization
$V_E^*(z):=\limsup_{\zeta \to z} V_E(\zeta)$ is called the global
extremal function of $E$; either $V_E^*\equiv +\infty$ -- this
occurs precisely when $E$ is {\it pluripolar}; i.e., $E\subset
\{u=-\infty\}$ for some $u\not \equiv -\infty$ psh on a neighborhood
of $E$ -- or else $V_E^*\in L^+(\C^N)$. It is well-known that if $E$
is a compact set in $\C^N$, then $V_E$ defined in (\ref{ve})
coincides with $V_E$ in formula (1) (cf., \cite{K} Theorem 5.1.7)
and hence $V_E$ is lowersemicontinuous. Thus for compact sets $E$,
$E$ is regular if and only if $V_E=V_E^*$.

As well as the classes $L(\C^N)$ and $L^+(\C^N)$, we will consider
the class
$$H(\C^N):=\{u\in L(\C^N): u(\lambda z)=u(z)+\log {|\lambda|} \ \hbox{for} \ \lambda \in \C, \ z\in \C^N\}$$
of {\it logarithmically homogeneous} psh functions.

Given $u:\C^N\to\R$ in $L(\C^N)$ we define the {\it Robin function}
of $u$ to be
$$\rho_u(z):=\limsup_{|\lambda|\to \infty} \left[u(\lambda z)-
\log |\lambda| \right].$$ Note that for $\lambda \in \C$,
$\rho_u(\lambda z)=\log {|\lambda|} + \rho_u(z)$; i.e., $\rho_u$ is
logarithmically homogeneous. It is known (\cite{bloom}, Proposition
2.1) that for $u\in L(\C^N)$, the Robin function $\rho_u(z)$ is
plurisubharmonic in ${\C}^N$; indeed, either $\rho_u\in H(\C^N)$ or
$\rho_u\equiv -\infty$. As an example, if $p$ is a polynomial of
degree $d$ so that $u(z):={1\over d}\log {|p(z)|}\in L(\C^N)$, then
$\rho_u(z)={1\over d}\log {|\hat p(z)|}$ where $\hat p$ is the top
degree $(d)$ homogeneous part of $p$. For a compact set $K$, we
denote by $\rho_K$ the Robin function of $V_K^*$; i.e., $\rho_K :=
\rho_{V_K^*}$.

Suppose now that $K$ is {\it circled}; i.e., $z\in K$ if and only if
$e^{it}z\in K$. Then the extremal function $V_K$ in (1) can be
gotten via
$$V_K(z)=\max [0,\sup \{u(z):u \in H(\C^N), \ u\leq 0 \ \hbox{on} \
    K\}]$$
$$=\max [0,\sup \{{1\over{\rm deg}p}\log {|p(z)|}: p \
\hbox{homogeneous polynomial}, \ ||p||_K\leq 1\}]$$ (\cite{K},
Theorem 5.1.6). Moreover, we have the following.

\begin{lemma} Let $K\subset \C^N$ be compact, circled, and nonpluripolar.
Then \beq\label{circled1} V_K^*(z)=\max [0,\rho_K(z)] \eeq and
\beq\label{circled2} {\rm supp}(dd^cV_K^*)^N\subset \{\rho_K=0\}.
\eeq
\end{lemma}

\begin{proof} Equation (\ref{circled1}) follows from the above equation for $V_K$, which shows that $V_K^*(\lambda z)=V_K^*(z)+\log {|\lambda|}$ provided $z,\lambda z \not \in \hat K$, and the definition of $\rho_K$: if $V_K^*(z) >0$, then
$$\rho_K(z):=\limsup_{|\lambda|\to \infty}[V_K^*(\lambda z)-\log {|\lambda|}]$$
$$=\limsup_{|\lambda|\to \infty}[V_K^*(z)+\log {|\lambda|}-\log {|\lambda|}]=V_K^*(z).$$
We have $\rho_K\in H(\C^N)$ and $\rho_K(z)=V_K^*(z)$ if $V_K^*(z)
>0$; since the set $\{z\in \C^N: \rho_K(z)\leq 0\}$ differs from
$\hat K=\{z\in \C^N: V_K(z)= 0\}$ by at most a pluripolar set,
(\ref{circled1}) follows (cf., Corollary 5.2.5 \cite{K}). The Robin
function $\rho_K$ is locally bounded away from the origin $0$ which
implies, by the logarithmic homogeneity, that $(dd^c\rho_K)^N=0$ on
$\C^N \setminus \{0\}$ (see section 4 for a discussion of the
complex Monge-Amp\`ere operator). This gives (\ref{circled2}).
\end{proof}
\bigskip

Let $u\in L(\C)$ and $d\mu(t)={i\over 4\pi}\Delta u(t)dt \wedge
d\bar t$  be its Riesz measure. Jensen's formula yields that
$\mu(\C):=\int_{\C} d\mu(t) \leq 1$. If, in addition, $u(0)=0$, we
have
$$u(t)=\int_{\C} \log |1- \frac t \zeta|
d\mu(\zeta)$$ (\cite{ronkin}, p. 37). In the notation introduced in
this section, Theorem 1.2 yields the following version of a
one-variable approximation result:

\begin{theorem} Let $u\in L^+(\C)\cap C(\C)$ with the additional property that
$$\lim_{|t|\to \infty} [u(t)- \log {|t|}]$$
exists. Given $\epsilon >0$, there exist polynomials $p_n,q_n$ of
degree $n=n(\epsilon)$ with \beq\label{app}u(t)-\epsilon \leq
\max\bigl[{1\over n} \log {|p_n(t)|}, {1\over n} \log
{|q_n(t)|}\bigr]\leq u(t), \ t\in \C. \eeq
\end{theorem}

Note that (\ref{app}) implies that $p_n$ and $q_n$ have no common
zeros; this will also follow from the proof of the theorem. This
immediately gives an approximation result for the class $H(\C^2)$ of
logarithmically homogeneous psh functions in $\C^2$.

\begin{corollary} Let $U\in H(\C^2)$ be logarithmically homogeneous with the additional property that
$u(t):=U(1,t)$ satisfies the hypotheses of the previous theorem.
Given $\epsilon >0$, there exist homogeneous polynomials $P_n,Q_n$
of degree $n=n(\epsilon)$ with no common factors such that
\beq\label{app2}U(z,w)-\epsilon \leq \max\bigl[{1\over n} \log
{|P_n(z,w)|}, {1\over n} \log {|Q_n(z,w)|}\bigr]\leq U(z,w), \
(z,w)\in \C^2. \eeq
\end{corollary}

\begin{proof} If (\ref{app}) holds, define
$$P_n(z,w):=z^np_n(w/z) \ \hbox{and} \ Q_n(z,w):= z^nq_n(w/z).$$
Note that if $p_n,q_n$ are of degree exactly $n$; i.e., if
$$p_n(t)=a_0+a_1t+\cdots +a_nt^n  \ \hbox{and} \ q_n(t)=b_0+b_1t+\cdots +b_nt^n$$
with $a_nb_n \not =0$, then $P_n(0,w)=a_nw^n$ and $Q_n(0,w) =
b_nw^n$. Otherwise, we may have $P_n(0,w)\equiv 0$ and/or
$Q_n(0,w)\equiv 0$. Then, since $U(1,w/z)+\log {|z|}=U (z,w)$ for
$z\not =0$, (\ref{app}) implies
$$U(z,w) -\epsilon \leq \max\bigl[{1\over n} \log {|P_n(z,w)|},
{1\over n} \log {|Q_n(z,w)|}\bigr]\leq U (z,w)$$ for $z\not =0$. But
$U$ is subharmonic on $z=0$ so
$$U(0,w)=\limsup_{z\to 0}U (z,w);$$
together with the previous inequalities, this yields (\ref{app2})
for all $(z,w)\in \C^2$.
\end{proof}

For a regular compact set $K\subset \C^N$, it is known that the
Robin function $\rho_K$ is continuous on $\C^N\setminus \{0\}$ (cf.,
\cite{bloom}). Thus, if $N=2$, $\rho_K(1,t)\in L^+(\C)\cap C(\C)$
and
$$\lim_{|t|\to \infty}[\rho_K(1,t)-\log |t|]=\lim_{|t|\to \infty}\rho_K(1/t,1)=\rho_K(0,1).$$
We can apply the corollary to $\rho_K$ to find homogeneous
polynomials $P_n,Q_n$ with \beq \rho_K(z,w)-\epsilon \leq
\max\bigl[{1\over n} \log {|P_n(z,w)|}, {1\over n} \log
{|Q_n(z,w)|}\bigr]\leq \rho_K(z,w).  \label{robinapprox} \eeq

To prove the first part of Theorem 1.1, for a regular circled set
$K\subset \C^2$, using (\ref{circled1}) from Lemma 2.1 and
(\ref{robinapprox}), we have \beq V_K(z,w)-\epsilon \leq
\max\bigl[{1\over n} \log {|P_n(z,w)|}, {1\over n} \log
{|Q_n(z,w)|},0\bigr]\leq V_K(z,w).  \label{extrapprox} \eeq This
gives uniform convergence of
$$\tilde u_n(z,w):=\max\bigl[{1\over n} \log {|P_n(z,w)|},
{1\over n} \log {|Q_n(z,w)|},0\bigr]\to V_K(z,w)$$ in Theorem 1.1.

For regular circled sets  $K\subset \C^2$, (\ref{circled2}) of Lemma
2.1 implies that
$${\rm supp}(dd^cV_K)^2\subset \{(z,w):\rho_K(z,w)=0\}.$$
We now show using (\ref{robinapprox}) and (\ref{extrapprox}) that
\beq U_n(z,w):=\max\bigl[{1\over n} \log {|P_n(z,w)-1|}, {1\over n}
\log {|Q_n(z,w)-1|}\bigr]\to V_K(z,w)  \label{betterapprox} \eeq
locally uniformly on $\C^2\setminus \{\rho_K=0\}$.

To prove (\ref{betterapprox}), we observe from the inequality
$|A-B|\leq 2\max [|A|,|B|]$ we have \beq\label{est}U_n(z,w) \leq
\max\bigl[{1\over n} \log {|P_n(z,w)|}, {1\over n} \log
{|Q_n(z,w)|},0\bigr] +{\log 2 \over n}. \eeq Now on a compact set
$E\subset \C^2\setminus \{\rho_K \leq 0\}$, by (\ref{extrapprox}),
given $\epsilon >0$ with $2\epsilon < \inf_E V_K$, for
$n>n_0(\epsilon)$,
$$\max \bigl[|P_n(z,w)|,|Q_n(z,w)|\bigr] > \exp {[n(V_K(z,w)-\epsilon)]} \ \hbox{on} \ E.$$
By choosing $n_0(\epsilon)$ larger, if necessary, we may assume
$$\exp {[n(V_K(z,w)-\epsilon)]} -1 > \exp {[n(V_K(z,w)-2\epsilon)]} \ \hbox{on} \ E$$
so that
$$\max \bigl[|P_n(z,w)-1|,|Q_n(z,w)-1|\bigr]   > \exp {[n(V_K(z,w)-2\epsilon)]} \ \hbox{on} \ E.$$
Together with (\ref{extrapprox}) and (\ref{est}), this proves local
uniform convergence outside of $ \{\rho_K \leq 0\}$. On compact
subsets of $ \{\rho_K<0\}$, the story is similar due to the
logarithmic homogeneity of $\rho_K$, ${1\over n} \log {|P_n(z,w)|}$,
and ${1\over n} \log {|Q_n(z,w)|}$ and  (\ref{robinapprox}): for
$r>0$, if $E:=\{z\in K:\rho_K(z) <-r\}$, by  (\ref{robinapprox}),
given $\epsilon >0$ with $\epsilon <r$, for $n>n_0(\epsilon)$,
$$\max \bigl[|P_n(z,w)|,|Q_n(z,w)|\bigr] < \exp {[-n(r-\epsilon)]} \ \hbox{on} \ E.$$
Thus, $|P_n(z,w)-1|,|Q_n(z,w)-1|> 1- \exp {[-n(r-\epsilon)]} \
\hbox{on} \ E$. We conclude that
$$\max\bigl[{1\over n} \log {|P_n(z,w)-1|},
{1\over n} \log {|Q_n(z,w)-1|}\bigr] > {1\over n} \log \bigl[1- \exp {[-n(r-\epsilon)]}\bigr] \ \hbox{on} \ E.$$ Hence $U_n\to 0$ uniformly on
$E$.

Note that since we assume that $K$ is polynomially convex and
circled, we have that
 \beq
\partial K = \{(z,w):\rho_K(z,w)=0\}.
\label{boundary}
 \eeq

Here is an illustrative example of the reduction scheme: let
$K=\{(z,w)\in \C^2: |z|^2+|w|^2\leq 1\}$ be the closed unit ball in
$\C^2$. Then $V_K(z,w)=\log^+{(|z|^2+|w|^2)^{1/2}}$ and
$\rho_K(z,w)=\log{(|z|^2+|w|^2)^{1/2}}$ so that $\rho_K(1,t)={1\over
2}\log {(1+|t|^2)}$. Note that the support of $\Delta \rho_K(1,t)$
is all of $\C$, but that
$$\int_{\C} |\log {|t|}| \Delta \rho_K(1,t)< +\infty.$$
Thus, Theorem 2.2 provides a uniform approximation of the strictly
subharmonic function ${1\over 2}\log {(1+|t|^2)}$ by a function of
the form $$\max\bigl[{1\over n} \log {|p_n(t)|}, {1\over n} \log
{|q_n(t)|}\bigr].$$

To summarize: using the results of this section, in order to compete
the proof of  the first part of Theorem 1.1, it remains to prove the
one-variable approximation result, Theorem 1.2.

\vskip 5mm

\section{\bf Main approximation result.}
    \label{sec:4}

\vskip 3mm

In this section, we prove Theorem 1.2. We work exclusively in the
complex
 plane $\C$ with variable $z$. To recall  the notation, given a
 measure
$\mu$ in  $\C$  with $\mu(\C)=1$, we consider its logarithmic
potential

 \beq
 V(z) = \int_{\C} \log |1- \frac z \zeta|
d\mu(\zeta). \label{potential}
 \eeq

\noindent We assume

\beq \lim_{z\to \infty} [V(z) - \log |z|] \ \hbox{exists},
 \label{limit}
\eeq

\beq \int_{\C} |\log |z|| d\mu(z) < \infty,
 \label{convergence}
\eeq

\noindent and that $V(z)$ is continuous in $\C$.

\bigskip

\noindent{\bf Claim 1:} \\
{\em For each $\eps>0$ there exist  a number $N$ and polynomials
$P(z)$ and $Q(z)$ of degree $N$ such that } \beq |V(z) - \frac 1 N
\max \{ \log|P(z)|, \log|Q(z)|\}  |<\eps, \ \ \ z\in \C.
\label{max_claim} \eeq \vspace{5mm} In order to prove this statement
we shall prove the following result:

\bigskip

\noindent{\bf Claim 2:} \\
{\em For each $\eps>0$ there exists  a number $N$, polynomials
$P(z)$ and $Q(z)$ of degree $N$, and sets $E,F\subset \C$, $E\cap F
=\emptyset$
 such that }

\[
|V(z) - \frac 1 N  \log|P(z)||<\eps, \ \ \ z\in \C\setminus E,
\]
\beq V(z) +\eps > \frac 1 N  \log|P(z)|,\ \ \ z\in  E,
\label{approx_claim1} \eeq and
\[
|V(z) - \frac 1 N  \log|Q(z)||<\eps, \ \ \ z\in \C\setminus F,
\]
\beq V(z) +\eps > \frac 1 N  \log|Q(z)|,\ \ \ z\in F.
\label{approx_claim2} \eeq

\subsection{Pattern of the proof}

 {\em Step 1}: It follows from
(\ref{limit}) and also from continuity of $V$ that $V$ is uniformly
continuous in $\C$. Convolving if need be with an appropriate bump
function one may assume that $\mu$ has the form \beq d\mu(z) =
a(z)d\sigma(z), \label{smoothness} \eeq \noindent where $\sigma$ is
Lebesgue measure and $a \geq 0$ is a smooth function in $\C$. It
follows from (\ref{convergence}) that
\[
a(z) \to 0 \quad \mbox{as} \quad z \to \infty.
\]
Define \beq A:=\max_{z\in \C} a(z). \label{maxa} \eeq

\medskip

{\em Step 2}:  We reduce the problem to the case when $\mu$ has
compact support. Given a number $R>0$ we let $Q_R$ denote the square
\[
Q_R=\{z=x+iy; |x|,|y| < R\}.
\]
Given $\eta>0$ we find  an integer $M$ and a number $R$ so that \beq
\int_{\C\setminus Q_R} |\log |\zeta|| d\mu (\zeta) < \eta,
\label{logtail} \eeq \beq \mu(\C\setminus Q_R)=1/M <\eta,
\label{tail} \eeq and \beq \max_{|z|>R/3} a(z) \leq \eta.
\label{smalla} \eeq Denote the logarithmic potential from the
portion of $\mu$ outside $Q_R$ by
\[
 V_\infty(z) := \int_{\C\setminus Q_R} \log |1- \frac z \zeta|
d\mu(\zeta).
\]
Finally, set \beq {1\over M}\log {r_\infty} := \int_{\C\setminus
Q_R} \log | \zeta| d\mu(\zeta). \label{rinfty}
 \eeq
\noindent Note that $r_{\infty} >R/\sqrt{2}$.

\begin{lemma}\label{Lemma1} Let
\beq w_\infty\in \C, \ \ \ |w_\infty|= 10 r_\infty. \label{winfty}
\eeq Then
\[\left |V_\infty (z) - \frac 1 M \log \left |
1-\frac z {w_\infty} \right|\right | \leq \ C_1 \ \eta, \ \ z\not\in
E_{w_\infty},
\]
and \beq
 \frac 1 M \log \left |
1-\frac z {w_\infty} \right| \leq V_\infty (z) +
 \ C_2 \ \eta, \ z\in E_{w_\infty},
\label{inequalityinfty} \eeq where $C_1,C_2$ are constants
independent of $w_{\infty}$ and
\[
E_{w_\infty}= \{z: |z-w_\infty|< \frac 1 {20} |w_\infty|\}.
\]
\end{lemma}
\medskip

\noindent{\bf Remarks  1.} It is clear that $E_{w_\infty}\cap
Q_R=\emptyset$
 and also
that it is possible to chose two different points $w'_\infty$ and
$w''_\infty$ satisfying (\ref{winfty}) so that
$E_{w'_\infty}\cap E_{w''_\infty} =\emptyset $. \\
\noindent{\bf 2.} The values of the constants in this lemma depend upon $A$. \\
\noindent{\bf 3.} We use the notation $a \prec b$ to mean $a\leq Cb$
with $C$ a
 constant independent of all parameters except perhaps $A$ and $a
 \asymp b$
 to mean $a \prec b$ and $b \prec a$.

\medskip

{\em Step 3.} Define
\[
 V_0(z) := \int_{ Q_R} \log |z - \zeta| d\mu(\zeta).
\]
Given Lemma 3.1, it remains to approximate $V_0$ by a function of
the form ${1\over N} \log |P_N(z)|$, where $P_N$ is a polynomial of
degree $N$. In order to construct this approximation we need a
special partition of $Q_R$. Existence of the desired partitions is
ensured by a lemma due to R. Yulmuhametov \cite{yul}. We state this
result in a form which is adjusted to our situation. Let $\hat{\mu}$
denote the restriction of $\mu$ to $Q_R$. We have
$\hat{\mu}(Q_R)=(M-1)/M$.
 Given an integer
$k$ we split $Q_R$ into $k(M-1)$ pieces each of measure $1/Mk$.

\medskip

\noindent{\bf Lemma A.} {\em Given an integer $k>0$, there exists a
covering of $Q_R$
\[
Q_R = \cup_{l=1}^{(k-1)M} Q^{(l)},
\]
 and a splitting of
$\hat{\mu}$,
\[
\hat{\mu}=\sum_{l=1}^{(k-1)M} \mu^{(l)},
\]
with the following properties:
\begin{itemize}
\item
Each  $Q^{(l)}$ is a rectangle with sides parallel to the coordinate
axes such that the ratio of longest to shortest side does not exceed
3;
\item
 each point in $Q_R$ belongs to at most four distinct rectangles $Q^{(l)}$;
 \item  $\mbox{supp} \ \mu^{(l)} \subset  Q^{(l)}$;
\item
\beq \mu^{(l)}(Q^{(l)})= \frac 1 {kM}. \label{mul} \eeq
\end{itemize}
}

Fix such a partition. We look for a polynomial
 $P_k$ of degree $N:=k(M-1)$ of the form
\[
P_k(z)=\prod_{l=1}^{N} (z-\zeta^{(l)}),
\]
where the choice of the points $\{\zl\}_1^{k(M-1)}\subset Q_R$ is
related to the partition.

Let $d(l):=\mbox{diam}(\ql)$. We then have $\mbox{Area}(\ql)\asymp
d(l)^2$. In choosing $\{\zl\}_1^{k(M-1)}$, we first observe that, by
(\ref{smoothness}) and (\ref{maxa}), $d(l)$ cannot be too small:
\[
d(l)\geq \frac 1 {3(MA)^{1/2}} \frac 1 {k^{1/2}}.
\]
We split the set of indices into two subsets: \beq
\mathcal{I}_k=\{l: 1\leq l \leq N, d(l)\leq
 k^{1/3}\frac 1 {3(MA)^{1/2}} \frac 1 {k^{1/2} }\}, \
\mathcal{J}_k=\{1,2, \ldots , N\}\setminus \mathcal{I}_k.
\label{normalsquares} \eeq We say that $\ql$ is a {\em normal}
rectangle  if $l\in \mathcal{I}_k$. For such rectangles we set \beq
\zl_0=kM \int_{\ql} \zeta d\ml(\zeta), \label{centersmass} \eeq the
center of mass of $\ml$ in $\ql$, and then take
\[
\zl:=\zl_0+\dl,
\]
where $\dl$ are any complex numbers satisfying \beq |\dl|\leq
k^{-5}. \label{deviations} \eeq

For $l\in \mathcal{J}_k$ we let $\zl\in \ql$ be any points of $Q_R$
 with the property that
\[
|\zl - \zeta^{(m)}|>k^{-5}, \ l,m\in \mathcal{J}_k, \ l\neq m.
\]

The choice of $\zl$'s is related to the integer $k$ and to the
corresponding partition; hence we write
\[
Z_k:=\{\zl\}_1^N, \quad E_k=\{z\in \C; \mbox{dist}(z,
Z_k)<k^{-10}\}.
\]

\medskip

{\em Step 4:} We approximate the finite potential $V_0$.

\begin{lemma}\label{Lemma2} For each $\eta> 0$ one can chose $k$ large enough
so that
\[
|V_0(z)- \frac 1 {N_k} \log|P_k(z)||<\eta, \ z\not\in E_k; \
V_0(z)+\eta > \frac 1 {N_k}  \log|P_k(z)|, \ z\in  E_k.
\]
\end{lemma}

Together with Lemma \ref{Lemma1} this statement immediately yields
 Claim 2 since
it allows us to chose two polynomials of the
 form $$(1-z/w'_\infty)^{N}[P_k(z)+C]^M, \
 (1-z/w''_\infty)^{N}[Q_k(z)+C]^M$$
such that the corresponding exceptional sets are disjoint.

We now give the proofs of lemmas \ref{Lemma1} and \ref{Lemma2}. We
begin with the  atomization of the external part of the potential,
 $V_{\infty}$; i.e., we prove Lemma \ref{Lemma1}.

\subsection{Proof of Lemma 3.1.}

The quantity to be estimated
\[
D_\infty (z)= V_\infty(z) - \frac 1 M \log \left | 1-\frac
z{w_\infty}\right |,
\]
admits two  representations: \beq D_\infty (z)= \int_{\C\setminus
Q_R}
  \left (
\log \left |1- \frac z \zeta \right |
   -\log \left |1- \frac z {w_\infty} \right |
  \right )d\mu (\zeta);
\label{inftyrepresentation1} \eeq and also
\[
D_\infty (z)   =
  \int_{\C\setminus Q_R} (\log |z-\zeta| - \log |z-w_\infty|)d\mu (\zeta) +
  \frac {\log 10} M =
\]
\beq \int_{\C\setminus Q_R} \log
    \left |
    1+\frac{w_\infty-\zeta}{z-w_\infty}
    \right |d\mu (\zeta) + \frac {\log 10} M.
\label{inftyrepresentation2} \eeq The term $ \frac {\log 10} M$ does
not exceed $\eta \log 10$ and does not
 influence our estimates. We consider the following cases:
\medskip

\underline{Case 1: $|z|\leq R/2$.}\\
\medskip

In this case it suffices to use the representation
(\ref{inftyrepresentation1}) and note that for $\zeta \not \in Q_R$,
\[
\log 1/2 \leq \log \left |1- \frac z \zeta \right |,
   \log \left |1- \frac z {w_\infty} \right |
\leq \log 3/2.
\]

\medskip
\underline{Case 2: $R/2 \leq |z|\leq 3 |w_\infty|$.}
\medskip

Note that the set $E_{w_{\infty}}$ is contained in this annulus. We
still use
 the representation (\ref{inftyrepresentation1}) and estimate each summand
independently. We have
\[
\int_{\C\setminus Q_R}
  \log \left |1- \frac z \zeta \right | d\mu (\zeta) =
\]
\[
\left ( \int_{\zeta \in \C\setminus Q_R, |\zeta|<4|w_\infty|}
    + \int_{ |\zeta|>4|w_\infty|}
\right )
  \log \left |1- \frac z \zeta \right | d\mu (\zeta)
= S_1(z)+S_2(z).
\]
We then have
\[
S_1(z)=\int_{\zeta \in \C\setminus Q_R, |\zeta|<4|w_\infty|}
\log|z-\zeta|d\mu (\zeta) - \int_{\zeta \in \C\setminus Q_R,
|\zeta|<4|w_\infty|} \log|\zeta| d\mu (\zeta)
\]
\[
=S_{11}(z)+S_{12}.
\]
Note that $S_{12}$ is independent of $z$; from (\ref{logtail}), $
|S_{12}|\asymp \eta $. In order to estimate $S_{11}(z)$ we mention
that according to (\ref{smoothness}) and (\ref{smalla})
\[
\int_{|z-\zeta|<1}\log|z-\zeta|d\mu (\zeta) \asymp  \ \eta;
\]
this is used for (\ref{inequalityinfty}). In the rest of the set $\{
{\zeta \in \C\setminus Q_R, |\zeta|<4|w_\infty|}\}$ we have
\[
0< \log |z-\zeta| < 10 \log r_\infty.
\]
Using (\ref{rinfty}) and (\ref{tail}) we have $ |S_{11}|\asymp \eta.
$

\medskip

When estimating $S_2(z)$ it suffices to observe that the integrand
is bounded and then apply (\ref{tail}).

\bigskip

\underline{Case 3: $ |z|\geq 3 |w_\infty|$.}
\medskip

We now use (\ref{inftyrepresentation2}). We have
\[
D_\infty(z)= \left ( \int_{\zeta\not\in Q_R, |\zeta| <2|w_\infty|}+
     \int_{2|w_\infty|<|\zeta|< 4|z|}+
         \int_{4|z|<|\zeta|}
             \right )
     \log
\left |
    1+\frac{w_\infty-\zeta}{z-w_\infty}
    \right |d\mu (\zeta)\]
\[+ \frac {\log 10} M
   = T_1(z)+T_2(z)+T_3(z)+ \frac {\log 10} M .
\]
We have $|T_1|\asymp 1/M $ since the integrand is bounded. When
estimating $T_2$ we observe that the integrand is bounded from above
throughout the whole region of integration thus it suffices to
estimate
 the integral
over the region $|z-\zeta|< |z|/5$, say, in which the integrand is
not bounded
 from below. In this domain we have
\[
\int_{|z-\zeta|< |z|/5}
     \log
\left |
    1+\frac{w_\infty-\zeta}{z-w_\infty}
    \right |d\mu(\zeta) =
\int_{|z-\zeta|< |z|/5} \log |z-\zeta|d\mu(\zeta) -
\]
\[
\int_{|z-\zeta|< |z|/5} \log |z-w_\infty|d\mu(\zeta).
\]
The estimate of the right hand side is similar to that of $S_1(z)$. Precisely, to get an {\it upper bound} on $\int_{|z-\zeta|< |z|/5} \log |z-w_\infty|d\mu(\zeta)$, since $ |z|\geq 3 |w_\infty|$ and $|z-\zeta|< |z|/5$, we have $|z-w_\infty|\leq 4|z|/3$ and $4|z|/5\leq |\zeta| \leq 6|z|/5$. Hence
$$\int_{|z-\zeta|< |z|/5} \log |z-w_\infty|d\mu(\zeta)\leq \int_{|z-\zeta|< |z|/5} \log (4|z|/3) d\mu(\zeta)$$
$$\leq  \int_{|z-\zeta|< |z|/5} \log (5|\zeta|/3) d\mu(\zeta).$$
From (\ref{logtail}) and (\ref{tail}), $\int_{|z-\zeta|< |z|/5} \log (5|\zeta|/3) d\mu(\zeta)|\asymp \eta $. For the other integral, 
$$\int_{|z-\zeta|< |z|/5} \log |z-\zeta|d\mu(\zeta)\asymp \eta$$
from (\ref{smoothness}) and (\ref{smalla}).

The estimate of $T_3$ is also straightforward; we use
$|z-w_\infty|>|w_\infty|$ and $|\zeta|>12|w_\infty|$ to obtain
\[
T_3(z)= \bigl |\int_{|\zeta|>4|z|}\log
\frac{|\zeta-w_\infty|}{|z-w_\infty|}d\mu (\zeta)\bigr |
    \leq
 \bigl |\int_{|\zeta|>4|z|}\log\frac{130|\zeta|}{12|w_\infty|}d\mu (\zeta)\bigr | \]
and apply (\ref{logtail}), (\ref{tail}), (\ref {rinfty}) and (\ref
{winfty}).

\subsection{Proof of Lemma 3.2.}

We turn to the atomization of the potential $V_0$.

We split the proof into several steps.

\noindent {\bf a.} Write
\[
D_0(z):=V_0(z)- \frac 1 {N} \log|P_k(z)|= \sum_{l=1}^N
\underbrace{\int_{\ql} (
             \log|z-\zeta| - \log |z-\zl|
                     )d\ml(\zeta)}_{j_l(z)}.
\]
We will estimate the contributions from $j_l$'s for $l\in
\mathcal{I}_k$ and $l\in \mathcal{J}_k$ separately. The general
estimate in {\bf b.} will be used in {\bf c}.

\medskip

\noindent {\bf b.} Estimation of $j_l(z)$: Assume $z\not\in \ql$.
Then
\[
j_l(z)= \Re \int_{\ql} (L(\zeta) -L(\zl))d\ml(\zeta)
\]
with
\[
L(\zeta)= \log(z-\zeta).
\]
Using the Taylor expansion
\[
L(\zeta)-L(\zl)=L'(\zl)(\zeta-\zl)+\int_{\zl}^\zeta
L''(s)(\zeta-s)ds=
\]
\[
L'(\zl)(\zeta-\zl_0)- L'(\zl)\dl+\int_{\zl}^\zeta L''(s)(\zeta-s)ds
\]
as well as  (\ref{centersmass}) and (\ref{mul}) we obtain
\[
j_l(z)= \frac {\dl}{Mk}(\frac 1 {z-\zl}) +
       \int_{\ql} \int_{\zl}^\zeta \frac {\zeta-s}{(z-s)^2}ds d\ml_\zeta.
\]
Taking (\ref{deviations}) into account we obtain \beq |j_l(z)|\leq
\frac 1 {Mk^{6}} \frac 1 {\ldist} + \frac 1 {kM} \frac
{d(l)^2}{\ldist^2}. \label{jestimate} \eeq

\medskip

\noindent {\bf c.} \underline{Contribution from remote normal rectangles.}\\
Consider \beq l\in \mathcal{I}_k \quad  \mbox{with} \quad \ldist >
3k^{-1/2}. \label{distant} \eeq It follows from the definition of
normal rectangle in (\ref{normalsquares}) and $l\in \mathcal{I}_k$
that
\[
|s-z|\prec k^{1/3}\ldist
\]
for all $s\in \ql$. Combining this with (\ref{jestimate}),
integrating with respect to  Lebesgue measure $\sigma$ over $\ql$,
and recalling that $\mbox{Area}(\ql)\asymp d(l)^2$, we obtain
\[
|j_l(z)|\prec \frac{k^{1/3}}{k^5}\int_{\ql}\frac{d\sigma(s)}{|s-z|}
+ \frac{k^{2/3}}k \int_{\ql} \frac{d\sigma(s)}{|s-z|^2}.
\]
\vspace{1cm} Therefore
\[
\sum_{l\in \mathcal{I}_k, \ \ldist > 3k^{-1/2}} |j_l(z)| \prec
k^{-14/3}\int_{|s-z|>3k^{-1/2}, |s|<2R}\frac {d\sigma(s)}{|s-z|}+
\]
\[
\frac{k^{2/3}}k \int_{|s-z|>3k^{-1/2},
|s|<2R}\frac{d\sigma(s)}{|s-z|^2} \prec k^{-1/3}\log k \to 0 \
\mbox{as} \ k\to \infty.
\]
Thus choosing $k$ large enough we can make the contribution from the
remote normal rectangles; i.e., those satisfying (\ref{distant}),
 arbitrarily small.

\medskip

\noindent {\bf d.}\underline{ Contribution from normal rectangles
  which are
 close to $z$.}\\
 Set
 \[
 \mathcal{B}_k(z):=\{ l\in \mathcal{I}_k:\ldist<3k^{-1/2}\}.
 \]
  In this section we  estimate
\[
\sum_{l\in \mathcal{B}_k(z)}j_l(z).
\]
It follows from the construction that the total number of indices in
$\mathcal{B}_k(z)$ is bounded by some constant independent of $z$
and $k$ and also, from the definition of normal rectangle, that all
 the rectangles $\ql$, $l\in \mathcal{B}_k(z)$ are contained
in the  disk $\{|\zeta-z|\leq C k^{-1/6}\}$, $C$ being independent
of $z$ and $k$. Let $\zeta^{(m)}$ be the point nearest to $z$ among
all $\{\zl\}_{l\in \mathcal{B}_k(z)}$. We then have, using
(\ref{smoothness}) and (\ref{maxa}),
\[
\sum_{l\in \mathcal{B}_k(z)}|j_l(z)|\prec
   \int_{\{|\zeta-z|\leq C k^{-1/6}\}}|\log|z-\zeta||d\sigma(\zeta)\]
   \[+
  | \log|z-\zeta^{(m)}|| \int_{\{|\zeta-z|\leq C k^{-1/6}\}}d\sigma(\zeta).
\]
Assuming now that $z\not\in E_k$ (i.e., $|z-\zeta^{(m)}|>k^{-10}$)
we obtain
\[
\sum_{l\in \mathcal{B}_k(z)}|j_l(z)|\prec k^{-1/3}\log k.
\]
Clearly if $z\in E_k$, we get a lower bound:
\[
\sum_{l\in \mathcal{B}_k(z)}j_l(z)\geq -C k^{-1/3}\log k.
\]
\medskip

\noindent {\bf e.} \underline{Contribution of non-normal rectangles.}\\
Define
\[
D_n(z):=\sum_{l\in \mathcal{J}_k}j_l(z).
\]
Let
\[
E=\cup_{l\in \mathcal{J}_k}\ql; \ \tilde{\mu}=\sum_{l\in
\mathcal{J}_k}\ml.
\]
From (\ref{normalsquares}), the area of each non-normal rectangle is
at least
 $(10MA)^{-1}k^{-1/3}$
and the total area they cover does not exceed $16R^2$ (since the
multiplicity of the covering is at most 4). Hence we have \beq
\sharp \mathcal{J}_k \prec k^{1/3}. \label{smallnumber} \eeq
Therefore
\[
\tilde{\mu}(Q_R) \prec k^{-2/3}.
\]
We first assume that $|z|<2R$. Letting $\zeta_m$ denote the point
which is the nearest to $z$ among all $\zl$, $l\in \mathcal{J}_k$,
we have
\[
|D_n(z)|\prec \int_{Q_R}|\log|z-\zeta||d\tilde{\mu}(\zeta)+
|\log|z-\zeta_m||\int_{Q_R}d\tilde{\mu}(\zeta)=A_1(z)+A_2(z).
\]

Now by (\ref{smoothness}) and (\ref{maxa}),
\[
|A_1(z)|\prec
     A \int_{|\zeta-z|<k^{-5}}|\log|z-\zeta||d\sigma(\zeta) \]
     \[ +
        \log k \int_{|\zeta-z|>k^{-5}, \zeta\in Q_R}d\tilde{\mu}(\zeta)\prec
  k^{-2/3}\log k.
\]
Assuming $z\not \in E_k$ (i.e., $|z-\zeta_m|>k^{-10}$) we have
\[
|A_2(z)|\prec \log k \tilde{\mu}(Q_R) \prec k^{-2/3}\log k.
\]
Otherwise we get a one-sided bound. These inequalities
 complete the estimate of $D_n$ in the case $|z|<2R$.

If $|z|>2R$ we simply have
\[
D_n(z)=\sum_{l\in \mathcal{J}_k} \int_{\ql} \left (
    \log \left |1-\frac \zeta z \right | -
                 \log \left |1-\frac {\zl} z \right |
\right ) d\tilde{\mu}(\zeta),
\]
and since the integrands are bounded we obtain
\[
|D_n(z)|\prec k^{-2/3}, \ |z|>2R.
\]
This inequality completes our estimates.

\vskip 5mm

\section{\bf Convergence of the Monge-Amp\`ere measures.}
    \label{sec:5}

\vskip 3mm

We return to $\C^2$ with variables $(z,w)$. We use the notation
$d=\partial +\bar \partial$ and $d^c =i(\bar \partial - \partial)$
where, for a $C^1$ function $u$,
$$\partial u:={\partial u\over \partial z}dz+{\partial u\over \partial w}dw, \ \bar \partial u:={\partial u\over \partial \bar z}d\bar z+ {\partial u\over \partial \bar w}d\bar w$$ so that $dd^c =2i \partial \bar \partial$. For a $C^2$ function $u$,
$$(dd^cu)^2=16 \bigl[ {\partial^2 u \over \partial z  \partial \bar z}{\partial^2 u \over \partial w  \partial \bar w}-{\partial^2 u \over \partial z \partial \bar w}{\partial^2 u \over \partial w  \partial \bar z}\bigr]{i\over 2}dz\wedge d\bar z\wedge {i\over 2}dw\wedge d\bar w$$
is, up to a positive constant, the determinant of the complex
Hessian of $u$ times the volume form on $\C^2$. Thus if $u$ is also
psh, $(dd^cu)^2$ is a positive measure which is absolutely
continuous with respect to Lebesgue measure. If $u$ is psh in an
open set $D$ and locally bounded there, or, more generally, if the
unbounded locus of $u$ is compactly contained in $D$, then
$(dd^cu)^2$ is a positive measure in $D$ (cf., \cite{BT1},
\cite{dem}). We discuss aspects of this last statement that we need.

A psh function $u$ in $D$ is an usc function $u$ in $D$ which is
subharmonic on components of $D\cap L$ for complex affine lines $L$.
In particular, $u$ is a locally integrable function in $D$ such that
\beq \label{ddc} dd^cu=2i\bigl[{\partial^2 u \over \partial z
\partial \bar z}dz\wedge d\bar z+{\partial^2 u \over \partial w
\partial \bar w}dw\wedge d\bar w+{\partial^2 u \over \partial z
\partial \bar w}dz\wedge d\bar w+ {\partial^2 u \over \partial \bar
z \partial w}d\bar z\wedge dw\bigr] \eeq is a positive $(1,1)$
current (dual to $(1,1)$ forms); i.e., a $(1,1)$ form with
distribution coefficients. Thus the derivatives in (\ref{ddc}) are
to be interpreted in the distribution sense. Here, a $(1,1)$ current
$T$ on a domain $D$ in $\C^2$ is positive if $T\wedge (i\beta \wedge
\bar \beta)$ is a positive distribution for all $(1,0)$ forms $\beta
= adz+bdw$ with $a,b\in C^{\infty}_0(D)$ (smooth functions having
compact support in $D$). Writing the action of a current $T$ on a
test form $\psi$ as $<T,\psi>$, this means that
$$<T, \phi (i\beta \wedge \bar \beta)>\geq 0 \ \hbox{for all} \ \phi \in C^{\infty}_0(D) \ \hbox{with} \ \phi \geq 0.$$
For a discussion of currents and the general definition of
positivity, we refer the reader to Klimek \cite{K}, section 3.3.

Following \cite{BT1}, we now define $(dd^cv)^2$ for a psh $v$ in $D$
if $v\in L^{\infty}_{loc}(D)$ using the fact that $dd^cv$ is a
positive $(1,1)$ current with measure coefficients. First note that
if $v$ were of class $C^2$, given $\phi\in C^{\infty}_0(D)$, we have
 $$\int_D \phi (dd^cv)^2=  -\int_D d\phi \wedge d^c v \wedge dd^cv $$
 $$= -\int_D dv \wedge d^c \phi \wedge dd^cv= \int_D v dd^c \phi \wedge dd^c v$$
since all boundary integrals vanish. The applications of Stokes'
theorem are justified if $v$ is smooth; for arbitrary psh $v$ in $D$
with $v\in  L^{\infty}_{loc}(D)$, these formal calculations serve as
motivation to {\it define} $(dd^cv)^2$ as a positive measure
(precisely, a positive current of bidegree $(2,2)$ and hence a
positive measure) via
$$<(dd^cv)^2,\phi>:=\int_D v dd^c \phi \wedge dd^c v.$$
This defines $(dd^cv)^2$ as a $(2,2)$ current (acting on $(0,0)$
forms; i.e., test functions) since $vdd^cv$ has measure
coefficients. We refer the reader to \cite{BT1} or \cite{K} (p. 113)
for the verification of positivity of $(dd^cv)^2$. Also, the use of
Stokes' theorem is valid and hence, for simplicity, we will write
$<(dd^cv)^2,\phi>$ as $\int_D \phi (dd^cv)^2$.

Despite the fact that $L^1_{loc}(D)$ might appear to be the natural
topology in which to study psh functions, work of Cegrell and Lelong
(cf., \cite{K} section 3.8) yields that on, e.g., a ball $D$, for
any psh function $v\in  L^{\infty}_{loc}(D)$, there always exists a
sequence of continuous psh functions $\{v_j\}$ with $v_j \to v$ in
$L^1_{loc}(D)$ but $(dd^cv_j)^2=0$ for all $j$. In the locally
bounded category, however, the complex Monge-Amp\`ere operator is
continuous under (a.e.) monotone limits (cf., Bedford-Taylor
\cite{BT2} or Sadullaev \cite{sad}). A simpler argument shows that
local uniform convergence of a sequence of locally bounded psh
functions $\{v_j\}$ to $v$ implies weak-* convergence $(dd^c v_j)^2
\to (dd^cv)^2$: in case $v_j,v$ are smooth, given $\phi\in
C^{\infty}_0(D)$,
 $$\int_D \phi (dd^cv_j)^2=  \int_D v_j dd^cv_j \wedge dd^c \phi$$
 $$=  \int_D v dd^cv_j \wedge dd^c \phi+ \int_D (v_j-v) dd^cv_j \wedge dd^c \phi.$$
 The first term tends to  $\int_D v dd^cv \wedge dd^c \phi=\int_D \phi (dd^cv)^2$ since $dd^cv_j \to dd^cv$ as positive $(1,1)$ currents; from the uniform convergence $v_j \to v$, the family $\{dd^cv_j\}$ is locally uniformly bounded (cf.,  \cite{sad}) so that the second term goes to zero. In particular, we obtain the following result.

\begin{proposition}
    \label{th:4.1} Let $K\subset \C^2$ be a regular, polynomially convex compact set. Suppose $\{u_n\}\subset L^+(\C^2)$ converges uniformly to $V_K$ on $\C^2$. Then
$$(dd^c u_n)^2 \to (dd^cV_K)^2$$
weak-* as measures in $\C^2$. Thus with $K, \ \{\tilde u_n\}$ as in
Theorem 1.1,
$$(dd^c \tilde u_n)^2 \to (dd^cV_K)^2.$$
\end{proposition}

The functions $\{U_n\}$ of Theorem 1.1 are {\it not} locally
bounded, but they are in the classical Sobolev space
$W^{1,2}_{loc}(\C^2)$. Following \cite{BT1} as before -- but
altering the final application of Stokes' theorem -- we note that if
$v\in  W^{1,2}_{loc}(D)$ for some domain $D$, and $\phi\in
C^{\infty}_0(D)$, we can formally write
 $$\int_D \phi (dd^cv)^2= - \int_D d\phi \wedge d^c v \wedge dd^cv $$
 $$=- \int_D dv \wedge d^c \phi \wedge dd^cv=- \int_D dv \wedge d^c v \wedge dd^c\phi$$
since all boundary integrals vanish. In this case, these
calculations serve as motivation to define $(dd^cv)^2$ as a positive
measure for a psh function $v$ in $W^{1,2}_{loc}(D)$ via
$$\int_D \phi (dd^cv)^2:=-\int_D dv \wedge d^c v \wedge dd^c\phi.$$
The functions $u(z,w):={1\over 2}\log {(|z|^2+|w|^2)}$ and $\tilde
u(z,w)=\max [\log {|z|},\log {|w|}]$ are canonical examples of such
functions with \beq\label{monge1} (dd^cu)^2=(dd^c\tilde u)^2=
(2\pi)^2\delta_{(0,0)} \eeq (\cite{dem}, Corollary 6.4). More
generally, if $f$ and $g$ are holomorphic functions near $(0,0)$, an
elementary calculation (cf., \cite{BT1}, p. 15) shows that
\beq\label{nomass}\bigl(dd^c {1\over 2} \log
{(|f|^2+|g|^2|)}\bigr)^2=0 \ \hbox{on} \ \{|f|^2 +|g|^2 > 0\}. \eeq
Thus if $f(0,0)=g(0,0)=0$ and $(0,0)$ is an isolated zero of
$\{f=g=0\}$, in a neighborhood of the origin, the Monge-Amp\`ere
measures
$$\bigl(dd^c \max (\log |f|,\log |g|)\bigr)^2, \ \bigl(dd^c {1\over 2} \log {(|f|^2+|g|^2|)}\bigr)^2$$
are supported at $(0,0)$. Indeed, we have
\beq\label{monge2}\bigl(dd^c \max (\log |f|,\log
|g|)\bigr)^2=\bigl(dd^c {1\over 2} \log
{(|f|^2+|g|^2|)}\bigr)^2=D(2\pi)^2\delta_{(0,0)} \eeq near $(0,0)$
where $D$ is the degree of the mapping $(z,w)\to (f(z,w),g(z,w))$ at
$(0,0)$. For example, taking $(z,w)\to (z,w^2)$,
$$\bigl(dd^c {1\over 2} \log {(|z|^2+|w|^4|)}\bigr)^2=2(2\pi)^2\delta_{(0,0)}.$$
To see how (\ref{monge1}) implies (\ref{monge2}), following
\cite{BT1}, p. 16, we observe that with $u(z,w):={1\over 2}\log
{(|z|^2+|w|^2)}$, the form
$$\omega:=  d^cu\wedge dd^cu$$
restricted to a sphere
$S_{\epsilon}:=\{(z,w):|z|^2+|w|^2=\epsilon^2\}$ equals
$2\epsilon^{-3}d\sigma_{\epsilon}$ where $d\sigma_{\epsilon}$ is the
volume form on $S_{\epsilon}$. If we write $F(z,w):=(f(z,w),g(z,w))$
and $v(z,w):={1\over 2} \log {(|f|^2+|g|^2|)}\bigr)^2$, then
$$d^cv\wedge dd^cv=F^*\omega = F^*(d^cu\wedge dd^cu).$$
Moreover,
$$\int F^*(\epsilon^{-3}d\sigma_{\epsilon})=2\pi^2D.$$
Hence
$$\int_{S_{\epsilon}}d^cv\wedge dd^cv=\int F^*(2\epsilon^{-3}d\sigma_{\epsilon})=4\pi^2D.$$
From (\ref{nomass}), $(dd^cv)^2$ is supported at $(0,0)$ and the
second equality in (\ref{monge2}) follows. The first follows from
Corollary 6.4 of \cite{dem}.

Thus for our functions
$$U_n(z,w)=\max\bigl[{1\over n} \log {|P_n(z,w)-1|},
{1\over n} \log {|Q_n(z,w)-1|}\bigr],$$ the Monge-Amp\`ere measures
$(dd^cU_n)^2$ are supported on the finite point sets
$K_n:=\{(z,w):P_n(z,w)=Q_n(z,w)=1\}$, and by the local uniform
convergence of $U_n \to V_K$ off of $\partial K=\{\rho_K=0\}$ (see
(\ref{boundary})), given $\epsilon >0$, for $n>n_0(\epsilon)$,
\beq\label{haus} K_n \subset (\partial
K)^{\epsilon}:=\{(z,w):|\rho_K(z,w)| \leq  \epsilon\}. \eeq From
Proposition 3.2 of \cite{blo}, in $\C^2$, convergence of a sequence
$\{v_j\}$ of psh functions in the Sobolev space
$W^{1,2}_{loc}(\C^2)$ implies weak-* convergence of the
Monge-Amp\`ere measures $\{(dd^cv_j)^2\}$; we will apply this result
to prove Theorem 4.1.

A simple example motivated from the one-variable example in the
introduction illustrates the distinction between approximation by
$\{\tilde u_n\}$ and by $\{U_n\}$. \vskip6pt

\noindent {\bf Example}. Let $K=\{(z,w):|z|,|w|\leq 1\}$ be the
closed unit bidisk. Then $$V_K(z,w)=\max [\log |z|, \log |w|,
0]=\max [\rho_K(z,w), 0]$$ so we can trivially take $P_n(z,w) =z^n$
and $Q_n(z,w)=w^n$ in Theorem 1.1. Then $\tilde u_n=V_K$ for all $n$
while
$$U_n(z,w)=\max [{1\over n}\log {|z^n-1|},{1\over n}\log {|w^n-1|}].$$
Thus $K_n$ consists of ordered pairs
$\zeta_{jk}^{(n)}:=(\omega_n^j,\omega_n^k), \ j,k=1,...,n$ where
$\omega_n =\exp {(2\pi i/n)}$ is a primitive $n-$th root of unity.
It is standard that
\begin{itemize}
\item $t\to t^n-1$ is a Fekete polynomial of degree $n$ for the closed unit disk in $\C$;
\item $(dd^c V_K)^2= d\theta_z \times d\theta_w$, the standard measure on the torus $T:=\{|z|=1\}\times \{|w|=1\}$ (of mass $(2\pi)^2$);
\item $U_n\to V_K$ locally uniformly in $\C^2\setminus K$ and $U_n \to 0$ locally uniformly in $K^o=\{\rho_k <0\}$, but $\{U_n\}$ does not converge pointwise on $T$; however,
\item $(dd^cU_n)^2={(2\pi)^2\over n^2}\sum_{j,k=1}^n\delta_{\zeta_{jk}^{(n)}}\to (dd^cV_K)^2$.
\end{itemize}
\vskip6pt

The assumption in Theorem 1.1 that $K$ is circled, regular and
polynomially convex implies that $K$ is {\it balanced}; i.e.,
$(z,w)\in K$ and $\lambda \in \C$ with $|\lambda| \leq 1$ imply
$(\lambda z,\lambda w)\in K$; moreover $K=\bar D$ where
$D=\{(z,w):\phi(z,w)<1\}$ is a balanced, pseudoconvex domain
determined by $\phi(z,w):=\exp {\rho_K(z,w)}$, the {\it Minkowski
functional} of $D$.

\begin{theorem} If $K=\bar D$ with $D$ strictly pseudoconvex, then
$$(dd^c U_n)^2 \to (dd^cV_K)^2$$
weak-* as measures in $\C^2$.
\end{theorem}

\begin{proof} We first note that all of the functions $U_n$ and $V_K$ have the same total Monge-Amp\`ere mass:
\beq\label{mamass}\int_{\C^2} (dd^c  U_n)^2 =  \int_{\C^2} (dd^c
V_K)^2 = (2\pi)^2. \eeq This is a standard fact about psh functions
$u\in L^+(\C^2)$; cf., \cite{taylor}.

Using \cite{H2}, Theorem 4.1.8, we can find a subsequence
$\{U_{n_j}\}$ of $\{U_n\}$ with $U_{n_j}\to U$ in $L^p_{loc}(\C^2)$
for some psh $U$ for all $p\in [1,\infty)$. Since $U_n \to V_K$
locally uniformly on $\C^2\setminus \{\rho_K=0\}$, $U_n\to V_K$ a.e.
in $\C^2$ and $\sup_n |U_n|$ is locally integrable. Hence $U=V_K$
and the full sequence $\{U_n\}$ converges; i.e.,  we have, in
particular, that $U_n \to V_K$ in both $L^2_{loc}(\C^2)$ and
$L^1_{loc}(\C^2)$. From this latter convergence, $\nabla U_n$
converges weakly (as distributions) to $\nabla V_K$. Using Blocki's
result, to show that $(dd^c U_n)^2 \to (dd^cV_K)^2$ weak-* as
measures, it thus suffices to show that $\nabla U_n \to \nabla V_K$
in $L^2_{loc}(\C^2)$. Note that $U_n,V_K\in W^{1,2}_{loc}(\C^2)$
(e.g., from \cite{blo}, Theorem 1.1).

Fix a strictly pseudoconvex domain $B=\{(z,w):\psi (z,w) <0\}$
containing $K$ where $\psi$ is strictly psh. We want to show that
$\nabla U_n \to \nabla V_K$ in $L^2(B)$. It suffices to show that
the norms converge; i.e.,
$$||\nabla U_n||^2:=\int_B |\nabla U_n|^2 \to \int_B |\nabla V_K|^2=||\nabla V_K||^2.$$
That is, by standard Hilbert space theory, weak convergence plus
convergence {\it of} the norms imply convergence {\it in} the norm.
Note that by the weak convergence of $\nabla U_n$ to $\nabla V_K$
(or simply Fatou's lemma) we have \beq\label{nab}\liminf_{n \to
\infty} ||\nabla U_n||\geq ||\nabla V_K||; \eeq we want to show the
limit exists and equals $||\nabla V_K||$.

Let $V_n:=\max[U_n,0]$. From the proof of the first part of Theorem
1.1 in section 2, $V_n\to V_K$ uniformly on $\C^2$ and hence, from
Proposition 1, $(dd^cV_n)^2 \to (dd^cV_K)^2$ weak-* as measures on
$\C^2$. By an observation of Cegrell, $V_n \to V_K$ in
$W^{1,2}_{loc}(\C^2)$. Precisely, {\sl If $\{u_j\}, \ u$ are
subharmonic functions in $W^{1,2}_{loc}(\R^m)$ and $u_j \to u$
locally uniformly, then $u_j \to u$ in $W^{1,2}_{loc}(\R^m)$}. To
see this, we may assume that $u_j,u$ are of class $C^2$ and we use
the identity
$${1\over 2} \Delta (v^2)=v\Delta v +|\nabla v|^2$$ for such functions. Take $\Omega'\subset \subset \Omega \subset \subset \R^m$ and $\eta \in C_0^{\infty}(\Omega)$ with $0\leq \eta \leq 1$ and $\eta =1$ on $\bar {\Omega'}$. Then
$$\int_{\Omega'}|\nabla (u_j-u)|^2 \leq \int_\Omega \eta |\nabla (u_j-u)|^2={1\over 2}\int_{\Omega} \eta \Delta [(u_j-u)^2]-\int_{\Omega} \eta  (u_j-u)\Delta (u_j-u)$$
$$\leq |{1\over 2}\int_{\Omega} (u_j-u)^2 \Delta \eta| +|\int_{\Omega} \eta (u_j-u)\Delta (u_j-u)|$$
$$\leq C\int_{\Omega} (u_j-u)^2+|\int_{\Omega} \eta (u_j-u)\Delta (u_j-u)|$$
(here $C$ depends on $\eta$) which tends to zero as $j\to \infty$
since $u_j \to u$ uniformly on $\bar {\Omega}$ and $\Delta u_j \to
\Delta u$ as measures.

We will work in an equivalent $L^2-$norm using a weight function. To
construct this function, we are assuming that $K=\bar D$ with
$D=\{(z,w):\rho_K(z,w) <0\}$ is strictly pseudoconvex; hence $\exp
{\rho_K}$ is strictly psh and we work on the sub-level sets
$B=B_R:=\{(z,w): \exp {\rho_K(z,w)}< e^R\}$ for $R>0$. For each set
$B$ we define $$\psi(z,w):=  \exp {\rho_K(z,w)}-e^R.$$ The (semi-)
norm in our new $L^2-$space is
$$||\nabla u||_{\psi}^2:= \int_B dd^c \psi \wedge d^cu \wedge du.$$
If $\psi(z) =A_1|z|^2 +A_2$ then $||\nabla u||^2_{\psi}=
4A_1||\nabla u||^2$; in general, due to strict plurisubharmonicity
and smoothness of $\psi$, we have
$$c_1 ||\nabla u|| \leq ||\nabla u||_{\psi} \leq c_2 ||\nabla u||$$
for constants $c_1,c_2$ depending only on $\psi$. The same argument
as before gives a version of (\ref{nab}) in our new norm:
\beq\label{nabwt}\liminf_{n \to \infty} ||\nabla U_n||_{\psi}\geq
||\nabla V_K||_{\psi}. \eeq

Now via integration by parts, we get
$$\int_B dd^c\psi \wedge dU_n \wedge d^cU_n =\int_B (-\psi) (dd^c U_n)^2 $$
modulo boundary integrals $\pm \int_{\partial B} dU_n \wedge d^cU_n
\wedge d^c\psi \pm \int_{\partial B} \psi d^cU_n \wedge dd^cU_n $.
Since $\psi =0$ on $\partial B$, this last term vanishes. Similarly,
$$\int_B dd^c\psi \wedge dV_K \wedge d^cV_K=\int_B (-\psi) (dd^c V_K)^2 $$
modulo boundary integrals $\pm \int_{\partial B} dV_K \wedge d^cV_K
\wedge d^c\psi \pm \int_{\partial B} \psi d^cV_K \wedge dd^cV_K $;
again, this latter term vanishes. Thus we must show that
\beq\label{bint}\int_{\partial B} dU_n \wedge d^cU_n \wedge d^c\psi
 \to \int_{\partial B} dV_K \wedge d^cV_K \wedge d^c\psi
\eeq and \beq\label{vint} \int_B (-\psi) (dd^c U_n)^2\to \int_B
(-\psi) (dd^c V_K)^2 . \eeq

Using (\ref{haus}), given $\epsilon >0$, for $n> n_0(\epsilon)$ we
have $(dd^c U_n)^2$ is supported in $(\partial K)^{\epsilon}$, and
$$1-2\epsilon -e^R \leq \psi(z,w) \leq 1+2\epsilon -e^R$$ on this set so that
$$(2\pi)^2 (1-2\epsilon-e^R) \leq \int_B \psi (dd^c U_n)^2 \leq (2\pi)^2(1+2\epsilon-e^R).$$
Since $(dd^c V_K)^2$ is supported on $\partial K$ and, from
(\ref{mamass}), the total Monge-Amp\`ere mass of $V_K$ is
$(2\pi)^2$,  we have $\int_B (-\psi) (dd^c V_K)^2=(2\pi)^2(e^R-1)$
so that
$$|\int_B (-\psi) (dd^c U_n)^2 - \int_B (-\psi) (dd^c V_K)^2|\leq (2\pi)^2 2\epsilon$$
for $n>n_0(\epsilon)$. This gives (\ref{vint}).

To prove (\ref{bint}), we observe that for any fixed $R>0$, for $n$
sufficiently large, $U_n=V_n$ on $\partial B =\partial B_R$. Thus we
may replace $U_n$ by $V_n$ in (\ref{bint}). Now $(dd^cV_n)^2 \to
(dd^cV_K)^2$ weak-* and the support of $(dd^cV_n)^2$ is compactly
contained in $B$ for $n$ large so
$$\int_B (-\psi) (dd^c V_n)^2 \to \int_B (-\psi) (dd^c V_K)^2.$$
Since $V_n \to V_K$ in $W^{1,2}_{loc}(\C^2)$,
$$\int_B dd^c\psi \wedge dV_n \wedge d^cV_n\to \int_B dd^c\psi \wedge dV_K \wedge d^cV_K.$$
Via the previously described integration by parts, (\ref{bint})
follows.

\end{proof}

\medskip

\noindent {\bf Remark 1}. If $K$ is not strictly pseudoconvex, if we
can find $\tilde K =\bar {\tilde D}$ balanced with $\tilde D$
strictly pseudoconvex and with supp$(dd^c V_K)^2\subset \tilde K$,
the same argument works using the function $\tilde \psi(z,w) = \exp
\rho_{\tilde K}(z,w) -e^R$. For example, for the bidisk $K$,
supp$(dd^c V_K)^2$ is the unit torus which is contained in the ball
$\tilde K=\{(z,w): |z|^2+|w|^2 \leq 2\}$.
\medskip

\noindent {\bf Remark 2}. Let $\Omega$ be a bounded hyperconvex
domain in ${\bf C}^N$; i.e., there exists a negative psh function
$\psi$ in $\Omega$ with $\{z\in \Omega: \psi (z) \leq -c\} \subset
\subset \Omega$ for all $c>0$. A bounded psh function $v$ belongs to
the class ${\cal E}_0(\Omega)$ if $\lim_{z'\to z}v(z')=0$ for all
$z\in \partial \Omega$ and $\int_{\Omega}(dd^cv)^N< +\infty$.
Finally, a psh function $v$ in $\Omega$ belongs to the class ${\cal
F}(\Omega)$ if there exists a sequence of functions $v_j\in {\cal
E}_0(\Omega)$ with $\sup_j \int_{\Omega}(dd^cv_j)^N < +\infty$ which
decreases to $v$ on $\Omega$. A recent result of Cegrell \cite{Ce}
states the following: for a sequence $\{u_n\}\subset {\cal
F}(\Omega)$, if $u_n\to u\in {\cal F}(\Omega)$ in
$L^1_{loc}(\Omega)$ and if there exists a strictly psh function
$v\in {\cal E}_0(\Omega)$ such that $\lim_{n\to \infty}
\int_{\Omega}v(dd^cu_n)^N=\int_{\Omega}v(dd^cu)^N$, then
$(dd^cu_n)^N$ converges weak-* to $(dd^cu)^N$. The sequence
$\{u_n\}$ must lie in ${\cal F}(\Omega)$ in order that certain
integration by parts formulae are valid. Note that functions in
${\cal E}_0(\Omega)$ have zero boundary values; moreover, if $u_n
\in {\cal F}(\Omega)$ then $\limsup_{z'\to z}u_n(z')=0$ for all
$z\in \partial \Omega$ (cf., \cite{ahag}). It might appear that (\ref{vint}) would
suffice (without (\ref{bint})) to prove Theorem 4.1. However, the
functions $U_n$ do not lie in the class ${\cal F}(B)$ since
$\limsup_{z'\to z}U_n(z')\not \equiv 0$ for all $z\in \partial B$.

\medskip

As mentioned in the introduction, from Bishop's construction, one
obtains the following result.

\begin{proposition} Let $K\subset \C^N$ be a regular, polynomially convex compact set. Then there exists a sequence of special polynomial polyhedra $\{\kappa_n\}$ where $\kappa_n$ is the closure of a union of a finite number of connected components of $${\cal K}_n:=\{(z_1,...,z_N): |P_{n,1}(z_1,...,z_N)| < 1, \ |P_{n,N}(z_1,...,z_N)| < 1\}$$
with $\{P_{n,1},...,P_{n,N}\}$ polynomials having degree $n$, such
that the extremal functions $\{V_{\kappa_n}\}$ converge uniformly to
$V_K$ and $(dd^c V_{\kappa_n})^N\to (dd^cV_K)^N$ weak-*.
\end{proposition}

However, it is not known how one can construct full component sets
of the form ${\cal K}_n$ approximating $K$ as we have in Theorem 1.1
using (\ref{calk}) nor how to construct functions $u_n$ of the form
$$u_n(z_1,...,z_N):= \max [{1\over n} \log {|\tilde
P_{n,1}(z_1,...,z_N)|}, ..., {1\over n} \log {|\tilde
P_{n,N}(z_1,...,z_N)|}]$$ for some polynomials $\tilde
P_{n,1},...,\tilde P_{n,N}$ so that, with
$$K_n:=\{(z_1,...,z_N):u_n(z_1,...,z_N)=-\infty\}$$ we have $(dd^c
u_n)^N$ is supported in $K_n$ as in (\ref{kd}) and
\begin{itemize}
\item $u_n \to V_K$ locally uniformly in $\C^N\setminus K$;
\item  $u_n \to V_K$ in $L^1_{loc}(\C^N)$; and
\item  $(dd^c u_n)^N\to (dd^cV_K)^N$ weak-*.
\end{itemize}

As a step in this direction, we can achieve a partial result in
$\C^2$.

\begin{proposition} Let $K\subset \C^2$ be a regular, polynomially convex compact set. Then there exists a sequence of pairs of polynomials $\{\tilde P_n, \tilde Q_n\}$
with deg$\tilde P_n=$deg$\tilde Q_n=n$ such that the functions
$$v_n(z,w):= \max [{1\over n} \log {|\tilde P_n(z,w)|}, {1\over n}
\log {|\tilde Q_n(z,w)|}]$$ converge to $V_K$ in
$L^1_{loc}(\C^2\setminus K)$ and $\rho_{v_n} \to \rho_K$ uniformly
on $\C^2$. In particular, if $K$ has empty interior (e.g., if
$K\subset \R^2$), then $v_n\to V_K$ in $L^1_{loc}(\C^2)$.
\end{proposition}

\begin{proof} Form the Robin function $\rho_K$ of $V_K$ (see section 2) and construct the regular, polynomially convex, circled set
$$K_{\rho}:=\{(z,w)\in \C^2: \rho_K(z,w)\leq 0\}.$$
Apply Theorem 1.1 to obtain a sequence of pairs $\{P_n,Q_n\}$ of
homogeneous polynomials such that if $\epsilon >0$ is given, then
$$\rho_K(z,w) -\epsilon \leq \max [{1\over n} \log {|P_n(z,w)|}, {1\over n} \log {|Q_n(z,w)|}] \leq \rho_K(z,w)$$
for all $(z,w)\in \C^2$ if $n > n(\epsilon)$. Construct
$$\tilde P_n ={\rm Tch}_K P_n, \  \tilde Q_n ={\rm Tch}_K Q _n$$
where, for a homogeneous polynomial $H_n$ of degree $n$,
$${\rm Tch}_K H_n := H_n + H_{n-1}$$
with deg$H_{n-1}\leq n-1$ and $||{\rm Tch}_K H_n||_K\leq ||H_n
+R_{n-1}||_K$ for all polynomials $R_{n-1}$ of degree at most $n-1$.
By Theorem 3.2 of \cite{bloom},
$$\limsup_{n \to \infty} ||\tilde P_n||_K^{1/n} \leq 1, \  \limsup_{n \to \infty} ||\tilde Q_n||_K^{1/n} \leq 1.$$
Thus, given $\epsilon >0$, for $n>n(\epsilon)$ we have
$$ \max [||\tilde P_n||_K,  ||\tilde Q_n||_K]\leq (1+\epsilon)^n$$
so that the the functions
$$v_n(z,w):= \max [{1\over n} \log {|\tilde P_n(z,w)|}, {1\over n} \log {|\tilde Q_n(z,w)|}]$$
satisfy
\begin{itemize}
\item $v_n\in L(\C^2)$;
\item given $\epsilon >0$, there exist $N=N(\epsilon)$ with $v_n \leq \epsilon$ on $K$ for $n > N(\epsilon)$;
\item $\rho_{v_n} \to \rho_K$ uniformly on $\C^2$.
\end{itemize}

\noindent This last item follows since
$$\rho_{v_n}=\max [{1\over n} \log {|P_n(z,w)|}, {1\over n} \log {|Q_n(z,w)|}] .$$
By Theorem 2.2 of \cite{bloom2}, we conclude that $v_n \to V_K$ in
$L^1_{loc}(\C^2\setminus K)$.
\end{proof}

\end{document}